\documentclass[12pt, reqno]{article}%{amsart}

\usepackage{latexsym, amsfonts, amsthm, amscd}
\usepackage{amssymb}

   \newenvironment{acknowledgments}{
       \section*{Acknowledgments}
    }{\par\addvspace{12pt}}

\def\isoS{\vbox{\baselineskip 0pt \lineskip .5pt
   \ialign{$\mathsurround=0pt  \scriptstyle \hfil ## \hfil $\crcr
       \sim \crcr = \crcr}}}

	\theoremstyle{plain}
	\newtheorem{prop}{Proposition}[section]
        \newtheorem{theorem}[prop]{Theorem}
	\newtheorem{thm}[prop]{Theorem}
	\newtheorem{lemma}[prop]{Lemma}

	\newtheorem{definition}[prop]{Definition}

	\theoremstyle{remark}

	\newcommand{\A}{\mathcal{A}}

	\newcommand{\ot}{\otimes}
	\newcommand{\op}{\oplus}

	\newcommand{\al}{\alpha}

	\newcommand{\om}{\omega}

	\newcommand{\beq}{\begin{eqnarray}}
	\newcommand{\eeq}{\end{eqnarray}}
	\newcommand{\beqs}{\begin{eqnarray*}}
	\newcommand{\eeqs}{\end{eqnarray*}}

	\newcommand{\lan}{\langle}
	\newcommand{\ran}{\rangle}

	\newcommand{\bthm}{\begin{theorem}}
	\newcommand{\ethm}{\end{theorem}}

	\newcommand{\ph}[1]{\phantom{#1}}

	\newcommand{\nn}{\nonumber}
	\newcommand{\na}{\nabla}
	\newcommand{\ct}{\mathbin{\vrule width1.5ex height.4pt\vrule height1.5ex} }

	\newcommand{\Ga}{\Gamma}
	\newcommand{\La}{\Lambda}
	\newcommand{\DE}{\widetilde{D}}
	\newcommand{\we}{\wedge}
	
	\newcommand{\Cinf}{C^\infty}

	\newcommand{\End}{\mbox{End}}

\newcommand{\ldiag}[1]%
       {\makebox[0cm]{${\scriptstyle#1}\downarrow\phantom{\scriptstyle#1}$}}
\newcommand{\ldiagup}[1]%
       {\makebox[0cm]{${\scriptstyle#1}\uparrow\phantom{\scriptstyle#1}$}}
\newcommand{\rdiag}[1]%
       {\makebox[0cm]{$\phantom{\scriptstyle#1}\downarrow{\scriptstyle#1}$}}
\newcommand{\sediagr}[1]%
       {\makebox[0cm]{$\phantom{\scriptstyle#1}\searrow{\scriptstyle#1}$}}
\newcommand{\nediagr}[1]%
       {\makebox[0cm]{$\phantom{\scriptstyle#1}\nearrow{\scriptstyle#1}$}}
\newcommand{\rdiagup}[1]%
       {\makebox[0cm]{$\phantom{\scriptstyle#1}\uparrow{\scriptstyle#1}$}}
\newcommand{\swdiag}[1]%
       {\makebox[0cm]{$\phantom{\scriptstyle#1}\swarrow{\scriptstyle#1}$}}
\newcommand{\sediag}[1]%
       {\makebox[0cm]{${\scriptstyle#1}\searrow\phantom{\scriptstyle#1}$}}
\newcommand{\nediag}[1]%
       {\makebox[0cm]{${\scriptstyle#1}\nearrow\phantom{\scriptstyle#1}$}}

\newcommand{\comdia}[9]{%
\begin{array}{ccc}
#1 & \stackrel{#2}{\longrightarrow} & #3 \\
\ldiag{#4} & #5 & \rdiag{#6} \\
#7 & \stackrel{#8}{\longrightarrow} & #9
\end{array}}

\begin{document}
\title{BV-generators and Lie algebroids} 

\author{ S\'ebastien Mich\'ea and Gleb Novitchkov
\footnote{Research partially supported by NSF grant DMS03-06665.}\\
\\
        Department of Mathematics\\
         Pennsylvania State University \\
         University Park, PA 16802, USA\\
\\
{\sf email: michea@math.psu.edu, gvn@math.psu.edu }}                                                         
\maketitle
\begin{abstract}
%We study the conditions that an operator,
%defined on the spaces of degree 0 and 1 of a Gerstenhaber or BV algebra,
%has to satisfy so that we can find an extension that generates the whole algebra.
%When applied 
%it gives %$\Z$-graded 
Let $\A=\op_i \A^i$ be a Gerstenhaber algebra
generated by $\A^0$ and $\A^1$.  Given a degree $-1$ operator $D$
on $\A^0\op\A^1$, we find the condition on $D$ that makes $\A$
a BV-algebra.  Subsequently, we apply it  to the Gerstenhaber 
or BV algebra associated to a Lie algebroid and obtain a global 
proof of the correspondence between BV-generators and flat connections.
\end{abstract}
\setcounter{section}{0}
\section{Introduction}
\setcounter{equation}{0}
Batalin-Vilkovisky (BV) algebras arose originally from the quantization of gauge field theories \cite{BV}.
In recent years, 
there has been a great deal of interest in these algebras in connection with various
subjects such as string theory and operads \cite{Bouwknegt, Getz, Kimura, Lian,
Penkava, Stasheff, Zwiebach}.

The correspondence between BV or Gerstenhaber algebras and various geometric structures
on a vector bundle has been studied by several authors \cite{Gerstenhaber, KosmannMagri, MackenzieXu}.  
The fact that Gerstenhaber algebras correspond Lie algebroids and strong differential Gestenhaber algebras
correspond to Lie bialgebroids was indicated in \cite{Kosmann}. 
In \cite{Xu} Xu established an explicit correspondence between Lie algebroids equipped with a flat connection
and BV-algebra structure on the space of their multi-sections.
In the particular case of multivector fields, the correspondence was 
found earlier by Koszul \cite{Koszul}. It has also been generalized by Huebschmann 
 to Lie-Rinehart algebras \cite{Huebschmann}. 

In this work we show a general result about Gerstenhaber and BV-algebras which,
when applied to the Gerstenhaber or BV-algebra associated to a Lie algebroid,
gives a new proof of Xu's result.
\bigskip

Let us give first the definitions of the main concepts used throughout this article.\\

A {\it Gerstenhaber algebra} is a triple $\big(\A=\oplus_{k\geq0} \A^i,\, \cdot,\, [,]\,\big)$
such that:
\begin{enumerate}
\item $\A$ %\A^0+\A^1+\A^2+\dots$  
 is a graded vector space,
\item The degree zero multiplication $\cdot$ endows $\A$ with a super-commutative 
associative algebra structure
\beqs
\A^i\cdot \A^j\subseteq\A^{i+j}.
\eeqs
Super-commutativity means that for each $a\in\A^{|a|}, \, b\in\A^{|b|}$,
\beqs
a\cdot b=(-1)^{|a|\cdot|b|}b\cdot a.
\eeqs

\item The degree $-1$ bracket $[,]$ endows $\A$ with a super-Lie algebra structure
\beqs
[\A^i, \A^j]\subseteq\A^{i+j-1}
\eeqs
satisfying the super-Leibniz identity: for every 
$a\in\A^{|a|}, b\in\A^{|b|},c\in\A^{|c|}$,
\beqs
[a,b\cdot c]=[a,b]\cdot c + (-1)^{(|a|-1)|b|}b\cdot [a,c].
\eeqs
The bracket $[,]$ also satisfies the super-Jacobi identity: for every 
$a\in\A^{|a|}, b\in\A^{|b|},c\in\A^{|c|}$,
\beqs
(-1)^{(|a|-1)(|c|-1)}[[a,b],c]+(-1)^{(|b|-1)(|a|-1)}[[b,c],a]+(-1)^{(|c|-1)(|b|-1)}[[c,a],b]=0.
\eeqs

\end{enumerate}

From now on all operators are assumed to be linear.

An operator $D$ of degree $-1$ is said to be a \emph{Gerstenhaber generator} if for
every $a\in\A^{|a|}$ and $b\in\A$,
\beqs
[a,b]=(-1)^{|a|}\Big(D(a\cdot b)-Da\cdot b -(-1)^{|a|}a\cdot Db\Big).\label{g-condition}
\eeqs

A Gerstenhaber algebra is called {\it exact}, if there is a Gerstenhaber generator $D$
satsfying $D^2=0$. 
An exact Gerstenhaber algebra is often referred to as {\it Batalin-Vilkovisky algebra}
(or {\it $BV$-algebra} for short) and a generating operator of vanishing square is called a $BV$-generator.

In the case that $\A^i=\Ga(\we^i A)$ and $\cdot=\we$, where $(A, a, [,]_A)$ is a Lie algebroid, 
the triple $\big(\A=\oplus_{k\geq0} \Ga(\we^i A),\we, [,]_A\,\big)$ is called 
\emph{the Gerstenhaber algebra of the Lie algebroid $(A, a, [,]_A)$}.

In this work we consider a Lie algebroid $(A, a, [,]_A)$ and the Gerstenhaber algebra 
$\big(\A=\oplus_{k\geq0} \Ga(\we^i A),\we, [,]_A\,\big)$ associated to it.  We suppose
there exists a degree $-1$ operator $D$ defined on $A^1=\Ga(A)$ and $A^0=\Cinf(M)$. 
The intent of this work is to study the necessary and sufficient conditions that make 
$\big(\A=\oplus_{k\geq0} \Ga(\we^i A),\we ,\, [,]_A\,\big)$ exact with generating
operator $D$.  

The paper is organized as follows.  

In Section $2$, we define an extension $\DE$ of any degree $-1$ operator $D$ defined on $\A^0\oplus\A^1$
using the bracket of the Gerstenhaber algebra $\big(\A=\oplus_{k\geq0} \A^i,\, \cdot,\, [,]\,\big)$,
and check that the extension is well defined.
Subsequently, we determine the conditions an 
operator $D$ has to satisfy for the extension to be a Gerstenhaber generator.
Then we investigate the condition imposed on $D$ to make its extension a BV-generator.

In Section $3$, we consider the Gerstenhaber algebra
$\big(\A=\oplus_{k\geq0} \Ga(\we^i A),\we, [,]_A\,\big)$ of a Lie algebroid $(A, a, [,]_A)$.
A generating operator was found in \cite{Xu} using local coordinates. Here we define first
a degree $-1$ operator on $A^0\oplus A^1=\Cinf(M)\oplus\Ga(A)$ without reference to local coordinates
and extend it to an operator $\DE$ defined on the whole algebra using the result of the second section.
We then recover Xu's conditions for % imposed on 
 an $A$-connection $\na$ on the line bundle 
$\we^n A$  that make the Gerstenhaber algebra $\big(\A=\oplus_{k\geq0} \Ga(\we^i A),\we, [,]\,\big)$ 
exact. Finally, we establish isomorphism between homology and cohomology spaces
\beq
H_k(A,\na_0)\cong H^{n-k}(A,\mathbb{R})
\eeq
for some particular flat $A$-connection $\na_0$ using the operator $\DE$.% in the case when there exists a non-vanishing%
%volume form on $A$.%

%%%%%%%%%%%%%%%%%%%%%%%%%%%%%%%%%%%%%%%%%%%%%%%%%%
%
%    SECTION: Gerstenhaber algebra generators
%
%%%
\section{Gerstenhaber and BV-algebra generators}
\setcounter{equation}{0}
%%%%%%%%%%%%%%%%%%%%%%%%%%%%%%%%%%%%%%%%%%%%%%%%%%
%
%    SUBSECTION: extension of degree -1 operators
%
%%%
\subsection{Extension of degree -1 operators}

In this section we consider a Gerstenhaber algebra $\A$ generated as super-algebra 
by $\A^0$ and $\A^1$.
Suppose we are given a degree $-1$ linear operator $D$ defined on $\A^0\op\A^1$,
then we define its extension $\DE$ to $\A$ by:
\beq
\DE(a\cdot b)&=&(-1)^{|a|}[a,b]+\DE a\cdot b +(-1)^{|a|}a\cdot \DE b,\label{def:extension of D}\\
&&\forall a\in\A^{|a|},\, b\in \A^{|b|}, \mbox{and }|a|+|b|>1;\nn\\
\DE&=&D, \mbox{ when restricted to  } \A^0\op\A^1.\label{def:extension of D_2}
\eeq
To make sure that the preceding formula unambigously  defines a degree $-1$
operator on $\A$, we have to check that for any $a\in \A^{|a|}, b\in \A^{|b|}, c\in \A^{|c|}$ 
we have :
\begin{itemize}
\label{comp1}\item[1)]{ $\DE(a\cdot b)=(-1)^{|a||b|}\DE(b\cdot a)$,} as $\A$ is 
a super-commutative :
\beq &&\DE(a\cdot b)\nonumber\\
&&= (-1)^{|a|}[a,b]+\DE a\cdot b + (-1)^{|a|}a\cdot \DE b\nonumber\\
&&=(-1)^{|a|}(-1)^{(|a|-1)(|b|-1)-1}[b,a]+(-1)^{(|a|-1)|b|}b\cdot \DE a 
+ (-1)^{|a|}(-1)^{|a|(|b|-1)}\DE b\cdot a\nonumber\\
&&=(-1){|a|+|a||b|-|b|-|a|}[b,a]+(-1)^{|a||b|+|b|}b\cdot\DE a + (-1)^{|a||b|}\DE b\cdot a\nonumber\\
&&=(-1)^{|a||b|}(-1)^{|b|}[b,a]+(-1)^{|a||b|}(-1)^{|b|}b\cdot\DE a+(-1)^{|a||b|}\DE b\cdot a\nonumber\\
&&=(-1)^{|a||b|}\DE(b\cdot a).\nn
\eeq
\label{comp2}\item[2)]{ $\DE((a\cdot b)\cdot c)=\DE(a\cdot (b\cdot c))$,} as $\A$ 
is an associative algebra :
\end{itemize}
\smallskip
\beq
\DE((a\cdot b)\cdot c)=(-1)^{|a|+|b|}[a\cdot b , c]+\DE(a\cdot b)\cdot c+(-1)^{|a|+|b|}(a\cdot b)\DE c\nn.
\eeq
Using the super-Leibniz property,
\beq &&(-1)^{|a|+|b|}[a\cdot b , c]=(-1)^{|a|+|b||c|}[a,c]\cdot b+(-1)^{|a|+|b|}a\cdot[b,c].\nn
\eeq
and 
\beq 
\DE(a\cdot b)\cdot c &=&(-1)^{|a|}[a,b]\cdot c + \DE a\cdot b\cdot c +(-1)^{|a|}
a\cdot\DE b\cdot c\nn
\eeq
so
\beq
\DE((a\cdot b)\cdot c)&=&(-1)^{|a|+|b||c|}[a,c]\cdot b+(-1)^{|a|+|b|}a\cdot[b,c]+(-1)^{|a|}[a,b]\cdot c\nn \\
&&+\DE a\cdot bc +(-1)^{|a|}a\cdot\DE b\cdot c+(-1)^{|a|+|b|}ab\cdot\DE c\label{dabc_lhs}.
\eeq
On the other hand,
\beq
\DE(a\cdot (b\cdot c))&=&(-1)^{|a|}[a,b\cdot c]+\DE a\cdot bc + (-1)^{|a|}a\cdot \DE(bc)\nn \\
&=&(-1)^{|a|} [a,b]\cdot c + (-1)^{|a|+|b||c|}[a,c]\cdot b+(-1)^{|a|+|b|}a\cdot[b,c]\nn\\
&&+\DE a\cdot bc+(-1)^{|a|}a\cdot\DE b\cdot c+(-1)^{|a|+|b|}ab\cdot \DE c\label{dabc_rhs}.
\eeq
Expressions (\ref{dabc_lhs}) and (\ref{dabc_rhs}) are identical, thus
$\DE((a\cdot b)\cdot c)=\DE(a\cdot (b\cdot c))$ for all $a\in \A^{|a|}, b\in \A^{|b|}, c\in \A^{|c|}$.
Therefore the operator $\DE$ is well defined, and, since we assumed that $\A$ is generated 
by $\A^0$ and $\A^1$ as a super-subalgebra, $\DE$ is defined everywhere on $\A$.
%\end{proof}

\subsection{The Gerstenhaber algebra generator condition}
Here we establish  the conditions that a degree $-1$ operator $D$ has to satisfy 
so that its extension $\DE$ is a Gerstenhaber generator.

The way we extended $D$ to $\A$ ensures that the Gerstenhaber generator (\ref{g-condition}) condition is 
satisfied by $\DE$ on $\A^2\op\A^3\op\cdots$. It remains to impose the conditions on $\A^0\op\A^1$.\\
On $\A^1=\A^0\cdot\A^1$, condition (\ref{g-condition}) reads
\beq
\DE(a^0 \cdot a^1)&=&[a^0, a^1]+D a^0\cdot b + a^0 \cdot D a^1\nn\\
&=&[a^0, a^1]+a^0 \cdot D a^1,\label{precomp_cond}
\eeq
since $D$ is of degree $-1$, and on $\A^0=\A^0\cdot\A^0$ condition (\ref{g-condition}) is always 
fulfilled.

Thus we have the following result
\begin{prop}A degree $-1$ operator $D$ defined on $\A^0\op\A^1$ admits an extension
to whole Gerstenhaber algebra $\A$ (generated by $\A^0$ and $\A^1$) that is a Gerstenhaber generator if and
only if it satisfies the condition
\beq
D(a^0\cdot a^1)=[a^0, a^1]+a^0 \cdot D a^1.  \label{comp_cond}
\eeq
Furthermore the extension is given by (\ref{def:extension of D}).
\end{prop}

%%%%%%%%%%%%%%%%%%%%%%%%%%%%%%%%%
%
%  SUBSECTION:  BV generator conditions
%
%%%
\subsection{BV-generator conditions}
To be a $BV$-generator, $\DE$ must also satisfy 
\beq 
\DE^2=0.
\eeq
The following is the necessary and sufficient condition for the extension $\DE$ defined
by (\ref{def:extension of D}), (\ref{def:extension of D_2}) to be a BV-generator.
%%%%%%%%%%%%%%%%%%%%%%%%%%%%%%%
%
%
%  D^2=0 <=> \DE^2=0.
%
%
%%%
\begin{prop}\label{thm:necessary and sufficient conditions for D^2=0 } $\DE^2=0$ if and only 
if
\beq
D[a,b]=[Da,b]+[a,Db], \,\,\forall a,b\in\A^1.  \label{thm:D^2=0 _c1}
\eeq
\end{prop}
\begin{proof}
Using a direct computation, we have
\beq
\DE^2(a\cdot b)=(-1)^{|a|}\Big\{\DE[a,b] - [\DE a, b] - (-1)^{|a|-1}[a,\DE b]\Big\}
+\DE^2a\cdot b + a\cdot\DE^2b.\,\,\,\,\label{d2}\\
\forall a\in\A^{|a|}, b\in \A^{|b|}.\nn
\eeq
Hence if $\DE^2=0$, then 
\beq
\DE[a,b] = [\DE a, b] + (-1)^{|a|-1}[a,\DE b], \label{DE2}
\eeq
which reduces to (\ref{thm:D^2=0 _c1}) when $a,b\in\A^1$.\\
Let us assume conversely that (\ref{thm:D^2=0 _c1}) holds, then using the following lemma
and the anti-symmetry of the bracket, a double induction on the degree of $a$ and $b$
shows that (\ref{DE2}) holds (namely $\DE$ is a derivation of the bracket).
So (\ref{d2}) reduces to
\beq
\DE^2(a\cdot b)=\DE^2a\cdot b + a\cdot\DE^2b,\,\,\,\,\label{d2}
\eeq
and by choosing $b$ in $\A^1$ an induction on the degree of $a$ proves that $\DE^2=0$.
\end{proof}
%%%%%%%%%%%%%%%%%%%%%%%%%%%%%%%%%%
%
%
%    Key Tool Lemma
%
%
%%%

\begin{lemma}\label{lem: key tool lemma}Let $a\in\A^{|a|}, b\in\A^{|b|}, c\in\A^{|c|}$. If 
\beq
\DE[a,b]&=&[\DE a, b] + (-1)^{|a|-1}[a,\DE b]\nn
\eeq
and
\beq
\DE[a,c]&=&[\DE a, c] + (-1)^{|a|-1}[a,\DE c],\nn
\eeq
then
\beq
\DE[a,bc]&=&[\DE a, bc] + (-1)^{|a|-1}[a,\DE(bc)].\label{DEa.bc}
\eeq
\end{lemma}
\begin{proof}
We simply compute separately each term of (\ref{DEa.bc}):
\beqs
\DE[a,bc]&=&\DE\Big([a,b]\cdot c+(-1)^{(|a|-1)|b|}b\cdot[a,c]\Big)\\
&=&(-1)^{|a|+|b|-1}[[a,b], c]+\DE[a,b]\cdot c + (-1)^{|a|+|b|-1}[a,b]\cdot\DE c\nn\\
&&+(-1)^{|a||b|}\big[b, [a,c]\big]+(-1)^{(|a|-1)|b|}\DE b\cdot [a,c] + (-1)^{|a||b|}b\cdot \DE[a,c].
\eeqs
Then 
\beqs
[\DE a, bc]&=&[\DE a, b]\cdot c + (-1)^{|a||b|}b\cdot[\DE a,c],
\eeqs
and
\beq
[a,\DE(bc)]&=&\Big[a, (-1)^{|b|}[b,c] + \DE b\cdot c + (-1)^{|b|}b\cdot\DE c\Big]\nn\\
&=&(-1)^{|b|}\big[a,[b,c]\big]+\big[a,\DE b\cdot c\big]+(-1)^{|b|}\big[a, b\cdot\DE c\big]\nn\\
&=&(-1)^{|b|}\big[a,[b,c]\big]+\big[a,\DE b\big] \cdot c + (-1)^{(|a|-1)(|b|-1)} \DE b\cdot [a,c]\nn\\
&&+(-1)^{|b|}\big[a, b\big] \cdot\DE c + (-1)^{|a||b|}b\cdot [a, \DE c].
\eeq
Now we put it all together:
\beq
&&\DE[a,bc]-[\DE a, bc] - (-1)^{|a|-1}[a,\DE(bc)]\nn\\
&=&\Big\{(-1)^{|a|+|b|-1}[[a,b], c]+(-1)^{|a||b|}\big[b, [a,c]\big]\Big\}\nn\\
&&-\Big\{[\DE a, b]\cdot c + (-1)^{|a||b|}b\cdot[\DE a,c]\Big\}\nn\\
&&-(-1)^{|a|-1}\Big\{(-1)^{|b|}\big[a,[b,c]\big]+\big[a,\DE b\big] \cdot c + (-1)^{|a||b|+|a|+|b|+1} \DE b\cdot [a,c]\nn\\
&&+(-1)^{|b|}\big[a, b\big] \cdot\DE c + (-1)^{|a||b|}b\cdot [a, \DE c]\Big\}\nn\\
&=&\Big\{(-1)^{|a|+|b|-1}[[a,b], c]+(-1)^{|a||b|}\big[b, [a,c]\big]
+(-1)^{|a|+|b|}\big[a,[b,c]\big]\Big\}\label{termz0}\\
&&+\Big\{(-1)^{|a||b|}b\cdot\DE[a,c]- (-1)^{|a||b|}b\cdot[\DE a,c]
-(-1)^{|a|-1}(-1)^{|a||b|}b\cdot [a, \DE c]\Big\}\nn.
\eeq
We can check using the super-Jacobi identity that (\ref{termz0}) is zero. We thus obtain
\beq
&&\DE[a,bc]-[\DE a, bc] - (-1)^{|a|-1}[a,\DE(bc)]\nn\\
&&=\Big\{\DE[a,b] - [\DE a, b] - (-1)^{|a|-1}\big[a,\DE b \big]\Big\}\cdot c\hskip2in\nn\\
&&\ph{+}+(-1)^{|a||b|}b\cdot\Big\{\DE[a,c]- [\DE a,c]-(-1)^{|a|-1}[a, \DE c]\Big\}\label{bs_thm:db}.
\eeq
And the right-hand side is zero using the hypotheses of the lemma.
\end{proof}

%%%%%%%%%%%%%%%%%%%%%%%%%%%%%%%%%%%%%%%%%%%%
%
%
%            Section 3
%   Application to Lie algebroids 
%
%
%%%

\section{Application to Lie algebroids}
\setcounter{equation}{0}

In this section we consider a  Lie algebroid $(A, M, a)$ of rank $n$ over a manifold $M$ and 
the associated Gerstenhaber algebra $\big(\A=\bigoplus_{i\geq0}\A^i,\, \we,\, [,]\,\big)$, 
where $\A^i=\Ga(\we^i A)$ and $[,]$ is the generalized Schouten braket.

If $E\to M$ a vector bundle,
an $A$-{\it connection} on $E$ is an $\mathbb{R}$-linear map $\na$:
\beqs
\Ga(A)\ot\Ga(E)&\to&\Ga(E),\\
     X\ot s    &\to& \na_X s,
\eeqs
satisfying axioms similar to those of the usual 
linear connection, that is for each $f\in\Cinf(M),\, X\in\Ga(A),\, s\in\Ga(E)$,
\beqs
&&\na_{fX}s=f\na_{X}s,\\
&&\na_{X}(fs)=(a(X)f)s+f\na_X s.
\eeqs
The {\it curvature} $R$ of an $A$-connection $\na$ is the element in
$\Ga(\we^2 A^*)\ot\End(E)$ defined by
\beq
R(X,Y)=\na_X\na_Y-\na_Y\na_X-\na_{[X,Y]}, \,\,\,\forall X, Y\in\Ga(A).
\eeq
An $A$-connection is {\it flat} if $R(X,Y)=0,\,\,\forall X, Y\in\Ga(A)$.
Next we establish the correspondence between $A$-connections and $BV$-generators.

%%%%%%%%%%%%%%%%%%%%%%%
%
%   Definition of D.
%
%
\subsection{Definition of the generating operator}
Given an $A$-{connection} $\na$ on the line vector bundle $\we^{n} A$ over $M$,
 we define a degree $-1$ operator $D$ on $\A^0\op\A^1=\Cinf(M)\op\Ga(A)$ by
\beq
Df&=&0,\,\,\,\,\,\forall f\in\Cinf(M)\label{def:DX 2},\\
(DX)\La &=& L_X\La - \na_X\La,\,\,\,\,\,\forall X\in\Ga(A)\label{def:DX};
\eeq
where $\La$ is any element of $\Ga(\we^{n} A)$ and $L_X$ is the Lie derivative with respect to $X$.\\
This definition of $D$  is independent of the choice
of $\La\in\Ga(\we^{n} A)$ as for any $\La'=f\La$, where $f$ is a smooth function,
we have:
\beq
L_X\La' -\na_X\La'&=&L_X(f\La) -\na_X(f\La)\nn\\
&=&\Big\{(a(X) f)\La + f L_X \La\Big\}- \Big\{(a(X) f)\La + f \na_X \La\Big\}\nn\\
&=&f\Big\{L_X \La - \na_X \La\Big\}\nn\\
&=&f(DX)\La\nn\\
&=&(DX)\La',
\eeq
as $DX\in \Cinf(M)$.\\
As seen in the previous section, a Gerstenhaber generator extension $\DE$ of $D$ exists if
and only if $D$ satisfies condition (\ref{comp_cond}).\\
In the case of the Gerstenhaber algebra of a Lie algebroid $A$ this condition reads:
\beq
D(fX)=[f,X]+fDX,\,\,\,\,\,\,\,\,\,\forall f\in\Cinf(M), X\in\Ga(A).\label{extension_condition}
\eeq
Now we check that it holds:
\beq
D(fX)\La&=&L_{fX}\La-\na_{fX}\La\nn\\
&=&\Big\{fL_X\La-X\we(d f\ct\La)\Big\}-f\na_X\La\nn\\
&=&f\Big\{L_X\La -  \na_X\La\Big\} - X\we(d f\ct\La)\nn\\
&=&fD(X)\La - X\we(d f\ct\La).
\eeq
Here the contraction (or interior product) $\ct$ is defined by 
$(\al\ct\La)(\beta)=\La(\al,\beta)$ for each $\al\in\Ga(\we^k A^*)$, 
$\beta\in\Ga(\we^{n-k} A^*)$, and $d$ is the usual Lie algebroid coboundary.\\
As
\beq
X\we(d f\ct\La)&=&\lan d f, X\ran\La\nn\\
&=&(a(X)f)\La,\nn
\eeq
and $[f,X]=-[X,f]=-a(X)f$, we obtain:
\beq
D(fX)\La=fDX\cdot\La+[f,X]\La.\nn
\eeq
As this equation is satisfied for any $\La$ in $\Ga(\we^{n} A)$, $D$ satisfies
 condition (\ref{extension_condition}). We can therefore extend $D$, in a unique way,
to a Gerstenhaber generator $\DE$ of the whole algebra $\A=\oplus_i \Ga(\we^i A)$
by using (\ref{def:extension of D}).  
%%%%%%%%%%%%%%%%%%%%%%%%%%%%%%%%
%
% Flat A-connections and BV-generators correspondence
%
%
%%%
\subsection{Flat A-connections and BV-generators correspondence}
The next theorem obtained by Ping Xu \cite{Xu} as a generalization to
any Lie algebroid of a Koszul's result for tangent bundle Lie algebroid \cite{Koszul}
is proved here in a coordinate-free framework.

\begin{thm}\label{A-connection <=> BV}Let $(A,M,a)$ be a Lie algebroid and
$\A=\op_{i}\Ga(\we^i A)$ be its associated Gerstenhaber algebra.
Then there exists a one-to-one correspondence between $A$-connections on
the line bundle $\we^n\A$ and linear operators $\DE$ generating the
Gerstenhaber algebra bracket.  Under this correspondence, flat connections 
correspond to $BV$-generators on $\A$.
\end{thm}
In the preceding section we associated a Gerstenhaber generator to any connection
on $\we^n A$. Here we first show in the following lemma that we can recover the connection
from its associated Gerstenhaber generator, hence proving the one-to-one correspondence.\\

%%%%%%%%%%%%%%%%%%%%%%%%%%%%%%%%%%
%
% V /\ =-X/\D(/\)
%  X
%
\begin{lemma}\label{lem:connection<->D}
Let $\na$ be an $A$-connection on $\we^n\A$ and $\DE$ its associated
Gerstenhaber generator.
Then for $X\in \Ga(A), \,\, \La\in \Ga(\we^n \La)$ we have:
\beq
\na_X \La = - X\we \DE\La.
\eeq
\end{lemma}
\begin{proof}
Since $\we^n A$ 
%$\Ga(\we^n A)$ 
is a line bundle, $X\we\La=0$ for any $\La\in \Ga(\we^n A)$ and $X\in \Ga(A)$.
Hence 
\beqs
L_X\La&=&[X,\La]\\
&=&-\Big\{\DE(X\we\La)-\DE X\we\La+X\we \DE\La\Big\}\\
&=&\DE X\cdot\La -X\we \DE\La,
\eeqs
Using the definition (\ref{def:DX}) of $DX=\DE X$, we obtain:
\beqs
\na_X\La&=&L_X\La -DX\cdot\La\\
&=&\Big\{DX\cdot\La -X\we \DE\La\Big\}-DX\cdot\La\\
&=&-X\we \DE\La.
\eeqs
\end{proof}

To prove the correspondence between BV-generators and flat connections we simply
show that the curvature $R(X,Y)$ of the connection $\na$ vanishes if and only if the square
$\DE^2$ of its associated Gerstenhaber generator vanishes also.
The following proposition shows that if $\DE^2=0$ then $R(X,Y)=0$, and 
that if $R(X,Y)=0$ then $\DE^2$ vanishes on $\Ga(\we^n A)$. 

%%%%%%%%%%%%%%%%%%%%%%%%%%%%%%%%%%%%%%%%5
%
%  R(X,Y)=-2 X/\Y/\D^2\La.
% 
%
%

\begin{prop}\label{prop: connection expr}
Let $\La\in\Ga(\we^n A)$. The extension $\DE^2$ is linked to the curvature 
$R$ of $\na$ by the following relation:
\beq
R(X,Y)\La&=&-X\we Y\we \DE^2\La,\label{lem_conn_expr1}\\
&&\forall X, Y\in \Ga(A).\nn
\eeq
\end{prop}
\begin{proof}
\beq
\na_X\na_Y\La&=&-X\we \DE(\na_Y\La)\nn\\
&=&-X\we \DE(-Y\we \DE\La)\nn\\
&=&X\we \DE(Y\we \DE\La)\nn\\
&=&X\we \Big\{-[Y, \DE\La]+DY\we \DE\La - Y\we \DE^2\La \Big\}\nn\\
&=&-X\we[Y, \DE\La]+DY\we X\we \DE\La -X\we Y\we \DE^2\La\nn\\
&=&-X\we[Y, \DE\La]+DY\we\big(-\na_X\La\big) - X\we Y\we \DE^2\La\nn\\
&=&-X\we[Y, \DE\La]-\big(L_Y-\na_Y\big)\na_X\La - X\we Y\we \DE^2\La\nn\\
&=&-X\we[Y, \DE\La]+\na_Y\na_X\La -L_Y\na_X\La - X\we Y\we \DE^2\La.\label{lem_naXY1}
\eeq
Similarly,
\beq
\na_Y\na_X\La&=-&Y\we[X, \DE\La]+\na_X\na_Y\La -L_X\na_Y\La-Y\we X\we \DE^2\La.\label{lem_naYX1}
\eeq
Now, (\ref{lem_naXY1})-(\ref{lem_naYX1}) gives:
\beq
\Big\{\na_X\na_Y-\na_Y\na_X\Big\}\La&=&Y\we[X, \DE\La]-X\we[Y, \DE\La] -\Big\{\na_X\na_Y-\na_Y\na_X\Big\}\La\nn\\
&&+\Big\{L_X\na_Y-L_Y\na_X\Big\}\La-2X\we Y\we \DE^2\La.
\eeq
Thus,
\beq
2\Big\{\na_X\na_Y-\na_Y\na_X\Big\}\La&=&Y\we[X, \DE\La]-X\we[Y, \DE\La]+\Big\{L_X\na_Y-L_Y\na_X\Big\}\La\nn\\
&&-2X\we Y\we \DE^2\La.
\eeq
By Lemma \ref{lem_constr_1} below,
\beq
\big(L_X\na_Y-L_Y\na_X\big)\La=2\na_{[X,Y]}-\Big\{Y\we [X, \DE\La]-X\we [Y, \DE\La]\Big\}\nn,
\eeq
therefore
\beq
2\Big(\na_X\na_Y-\na_Y\na_X\Big)\La=2\na_{[X,Y]}\La-2X\we Y\we \DE^2\La\nn,
\eeq
or
\beq
\Big(\na_X\na_Y-\na_Y\na_X-\na_{[X,Y]}\Big)\La=-X\we Y\we \DE^2\La.
\eeq
This completes the proof.
\end{proof}
%%%%%%%%%%%%%%%%%%%%%%%%%%%%%%%%%%%%%%%%
%
%  (L_X\na_Y-L_Y\na_X)\la = 2\na_{[X,Y]} - (Y/\[X, D\La]-X/\[Y, D\La])
%
%

\begin{lemma}\label{lem_constr_1}
For any $X, Y\in \Ga(A)$ and $\La\in\Ga(\we^n A)$, we have the following identity
\beq
\big(L_X\na_Y-L_Y\na_X\big)\La=2\na_{[X,Y]}-\Big\{Y\we [X, \DE\La]-X\we [Y, \DE\La]\Big\}.
\eeq
\end{lemma}
\begin{proof}
\beq
\big(L_X\na_Y-L_Y\na_X\big)\La&=&[X,\na_Y\La]-[Y,\na_X\La]\nn\\
&=&[X,-Y\we \DE\La] - [Y,-X\we \DE\La]\nn\\
&=&[Y,X\we \DE\La] - [X,Y\we \DE\La]\nn\\
&=&[Y,X]\we \DE\La+X\we [Y, \DE\La]\nn\\
&&-\Big\{[X,Y]\we \DE\La+Y\we [X, \DE\La]\Big\}\nn\\
&=&-2[X,Y]\we \DE\La-\Big\{Y\we [X, \DE\La]-X\we [Y, \DE\La]\Big\}\nn\\
&=&2\na_{[X,Y]}-\Big\{Y\we [X, \DE\La]-X\we [Y, \DE\La]\Big\}.
\eeq
\end{proof}

The last step in the proof of the theorem is to show that if $\DE^2=0$ on $\Ga(\we^n A)$
then it vanishes on the whole Gerstenhaber algebra $\A$.

%%%%%%%%%%%%%%%%%%%%%%%%%%%%%%%%%%%%%%%%%%%55
%
% Submain theorem
%  D^2=0 and \DE^2\La=0 => \DE^2=0.
% 
%

\begin{prop}\label{prop: D^2=0 everywhere}
Suppose $\DE^2\La=0, \,\,\forall \La\in \Ga(\we^n A)$. Then
\beq
\DE^2 U=0, \,\,\forall U\in \A.
\eeq
\end{prop}
\begin{proof}
We first show that for all $X\in\A^1, \La\in \Ga(\we^n A)$,
\beq
\DE[X,\La]=[\DE X, \La]+[X,\DE\La]. \label{thm:DE2}
\eeq
Since $A$ is of rank $n$, we have
$X\we\La=0$. Therefore by (\ref{d2}) we have
\beqs
\DE^2(X\we\La)&=&-\Big\{\DE[X,\La]-[DX,\La]-[X,\DE \La]\Big\}+D^2X\we\La+X\we\DE^2\La=0.\nn
\eeqs
$D^2X=0$ because $D$ is of degree $-1$ and $\DE^2\La=0$ by assumption, so
\beqs
\DE[X,\La]-[DX,\La]-[X,\DE \La]=0.
\eeqs
Now for any $X, Y\in \Ga(A)$ and $\La\in\Ga(\we^n A)$ we also have $Y\we\La=0$, so
using formula (\ref{bs_thm:db}) we obtain
\beq
0&=&\DE[X,Y\we\La]-[\DE X, Y\we\La]-[X,\DE(Y\we\La)]\nn\\
&=&\Big\{\DE[X,Y] - [\DE X, Y] - \big[X,\DE Y \big]\Big\}\we\La%\hskip2in\nn\\
\ph{+}+(-1)^{1\cdot 1}Y\we\Big\{\DE[X,\La]- [D X,\La]-[X, \DE\La]\Big\}\nn\\
&=&\Big\{D[X,Y] - [D X, Y] - [X,DY]\Big\}\La.\label{thm:DE2_2}
\eeq
As (\ref{thm:DE2_2}) is zero for all $\La\in\Ga(\we^n A)$, we get
\beq
D[X,Y] - [D X, Y] - [X,DY]=0.\label{exp: nec and suf1}
\eeq
%Therefore (\ref{prop:DE2_1}) holds. 
Hence $\DE$ is a derivation of the bracket.
According to Proposition \ref{thm:necessary and sufficient conditions for D^2=0 },
this implies that $\DE^2$ vanishes identically.  
\end{proof}

%%%%%%%%%%%%%%%%%%%%%%%%%%%%%%%%%%%%%%%%%%%%%%%%%%%%%%%%%%%%%%%
%
%
%         Lie Algebroid Homology.
%
%
%

\subsection{Lie Algebroid Homology}

Let $(A, a, [,]_A)$ be a Lie algebroid of rank $n$ and $\na$ a flat connection on the line bundle $\we^n A$.\\
Let $\DE$ be the corresponding Gerstenhaber generator and $\partial=(-1)^{n-k}\DE$ when 
restricted to $\A^k=\Ga(\we^k A)$. 
As $\na$ is flat, $\partial^2=0$ and we get a chain complex.  We denote by $H_*(A,\na)$ 
its homology:
\beq
H_*(A,\na)=ker\,\partial/Im\,\partial.
\eeq

We establish a relation between the Lie algebroid homology $H_*(A,\na)$
and the Lie algebroid cohomology with trivial coefficients $H^{*}(A,\mathbb{R})$
 in the case where the line bundle $\we^n A$ is trivial.

\begin{definition}Let $\La\in\Ga(\we^n A)$ be a nowhere vanishing section.
We define the operator $*$ from $\Ga(\we^{n-k}A^*)$ to $\Ga(\we^k A)$ by
\beq
*\om&=&\om\ct\La,\\
\forall\ \om&\in&\Ga(\we^{n-k}\A^*).\nn
\eeq
\end{definition}
By assuming that $\we^n A$ is a trivial line bundle, {\it i.e.}, 
that there exists a nowhere vanishing section $\La\in\Ga(\we^n A)$, %then 
we can construct
 a flat $A$-connection $\na_0$ on $\we^n A$ by $(\na_0)_X\La=0$
for all $X\in \Ga(A)$.\\
The   $*$-operator defined above becomes an intertwiner between the homological 
and cohomological spaces:
\begin{thm}\label{Hom/Cohom pairing}
Let $\na_0$ be a flat $A$-connection on $\we^n A$ and $\DE_0$ be the associated
generating operator. Then 
\beq
\DE_0*\om=-(-1)^{n-k}*d\om,\,\,\,\,\,\forall \om\in \Ga(\we^{n-k}A^*),
\eeq
and the following diagram is commutative:
\beq
\comdia{\Ga(\we^{n-k}A^*)}{*}{\Ga(\we^k A)}{d}{}{-\partial_0}{\Ga(\we^{n-k+1} A^*)}{*}{\Ga(\we^{k-1} A)},
\label{Hom/Cohom pairing:diag}
\eeq
where $d$ is the usual Lie algebroid coboundary.
\end{thm}
This leads directly to a global proof of Theorem 4.6 of \cite{Xu} :
\begin{thm} Let $\na_0$ be an $A$-connection on $\we^n A$ that admits a global nowhere vanishing 
horizontal section $\La\in \Ga(\we^n A)$. Then
\beqs
H_*(A,\na_0)\isoS H^{n-*}(A,\mathbb{R}).
\eeqs
\end{thm}
To prove the theorem, we need 
the following statements that are easily proven by induction.
\begin{lemma}Let $X_1,\dots, X_k \in \Ga (A)$. Then for any $X\in\Ga (A)$ and $\al\in\Ga(A^*)$,
we have :
\beq
[X,X_1\we\dots\we X_k]=\sum_{i=1}^k (-1)^{i-1} [X,X_i]\we X_1\we\dots\we\hat{X_i}\we\dots\we X_k.
\eeq
\beq
\al\ct X_1\we\dots\we X_k=\sum_{i=1}^k (-1)^{i-1}(\al\ct X_i)X_1\we\dots\hat{X_i}\dots\we X_k.
\eeq
\end{lemma}

\begin{lemma}\label{al_jU} Let $X_1,\dots, X_n \in \Ga(A)$ be elements of a basis in $\Ga(A)$, and let 
$\al_j=X^*_j$ be a dual basis. 
Then the following holds:   for $1\leq j\leq k$,
\beq
\al_j\ct X_1\we\dots\we X_k=(-1)^{j-1}X_1\we\dots\hat{X_j}\dots\we X_k.
\eeq
For $j>k$, \beq\al_j\ct X_1\we\dots\we X_k=0\nn.\eeq
\end{lemma}

%%%%%%%%%%%%%%%%%%%%%%%%%%%%%%%%%%%%%%%%%%%%%%
%
% dw_|U=w_|DU+(-1)^w D(w_|U)
%
%
Since the definition of the algebroid coboundary is given locally,
the following general proposition, which relates the Lie algebroid coboundary
to its Gerstenhaber generator, is proved using local coordinates. Its result, however, is global.

\begin{prop}\label{prop:dw general}Let $(A, a, [,]_A)$ be a Lie algebroid and let $\DE$ 
be a generating operator of the Gerstenhaber algebra $\A=\sum_{k=0}^n \Ga(\we^k A)$. 
Then, for any section $U\in\Ga(\we^u A)$ and any $\om\in\Ga(\we^{|w|}A)$ with $|w|+1\leq u$, 
\beq
d\om\ct U=\om\ct \DE U -(-1)^{|\om|}\DE (\om\ct U).\label{prop:dw::main relation general.}
\eeq
\end{prop}
\begin{proof}
In the first step we prove the theorem when $|w|+1 = u$ and $U$ is a nowhere vanishing 
section.\\
Let us choose a basis $X_1,\dots, X_n$ of 1-sections such that $U=X_1\we\dots\we X_{k+1}$ 
and again denote its dual basis by
$\al_1,\dots, \al_n$.
For any $\om\in \Ga(\we^k A^*)$, the Lie algebroid coboundary
is defined (see \cite{MackenzieXu, Xu}) by
\beq
d\om\ct U&=&d\om(X_1, \dots, X_{k+1})=\sum_{i=1}^{k+1} (-1)^{i+1}a(X_i)(\om(X_1, .\,\hat{X_i}\,., X_{k+1}))\nn\\ 
&&\phantom{\sum}+\sum_{i<j}(-1)^{i+j}\om([X_i, X_j], X_1, .\,\hat{X_i}\,.\,.\,\hat{X_j}\,.\,, X_{k+1}).
\label{def_dw}
\eeq
Using Lemma \ref{al_jU}
\beq
(-1)^{i+1}a(X_i)(\om(X_1, .\,\hat{X_i}\,., X_{k+1}))&=&a(X_i)(\om(\al_i\ct U))\nn\\
&=&a(X_i)(\om\ct(\al_i\ct U))\nn\\
&=&[X_i,\om\ct(\al_i\ct U)].\nn
\eeq
Thus we can write
\beq
\sum_{i=1}^{k+1} (-1)^{i+1}a(X_i)(\om(X_1, .\,\hat{X_i}\,., X_{k+1}))
&=&\sum_{i=1}^{k+1}[X_i,\om\ct(\al_i\ct U)].\nn%\\
%&=&\sum_{i=1}^{n}[X_i,\om\ct(\al_i\ct U)].\nn
\eeq
Using the properties of contraction,
\beq
[X_i,\om\ct(\al_i\ct U)]&=&(-1)^{|w|}[X_i,\al_i\ct(\om\ct U)].\nn
\eeq
Next we express the bracket using the definition of the generating operator $\DE$:
\beq
(-1)^{|w|}[X_i,\al_i\ct(\om\ct U)]&=&(-1)^{|w|}(-1)
\Big\{\DE(X_i\we(\al_i\ct(\om\ct U))-\DE X_i\we(\al_i\ct(\om\ct U))\nn\\
&&+X_i\we \DE(\al_i\ct(\om\ct U)) \Big\}.\nn
\eeq
As $\al_i\ct(\om\ct U)$ is a function, $\DE(\al_i\ct(\om\ct U))=0$.  Therefore we get:
\beqs
(-1)^{|w|}[X_i,\al_i\ct(\om\ct U)]&=&
-(-1)^{|w|}\Big\{\DE\big(X_i\we(\al_i\ct(\om\ct U))\big)%\\
-\DE X_i\we(\al_i\ct(\om\ct U))\Big\}.
\eeqs
Since $\al_i$ is dual to $X_i$,
\beq
\sum_{i=1}^{k+1}X_i\we(\al_i\ct(\om\ct U))&=&\om\ct U.\nn
\eeq
This allows us to simplify the first term of (\ref{def_dw}):
\beq
&&\sum_{i=1}^{k+1} (-1)^{i+1}a(X_i)(\om(X_1, .\,\hat{X_i}\,., X_{k+1}))\nn\\
&=&\sum_{i=1}^{k+1}-(-1)^{|w|}\Big\{\DE\Big(X_i\we(\al_i\ct(\om\ct U))\Big)-\DE X_i\cdot(\al_i\ct(\om\ct U))\Big\}\nn\\
&=&-(-1)^{|w|}\DE (\om\ct U)+(-1)^{|w|}\sum_{i=1}^{k+1}\DE X_i\cdot(\al_i\ct(\om\ct U)).
\eeq
For the second term of (\ref{def_dw}), we can write
\beq
&&\sum_{i=1}^{k}\sum_{i<j}(-1)^{i+j}\om([X_i, X_j], X_1, .\,\hat{X_i}\,.\,.\,\hat{X_j}\,.\,, X_{k+1})\nn\\
&=&\om\ct\sum_{i=1}^{k}(-1)^i X_1\we \dots\we X_{i-1}\we[X_i,V_i],\nn
\eeq
where $V_i=X_{i+1}\we \dots \we X_{k+1}$.\\
As $X_i\we V_i=V_{i-1}$ and $\DE$ is a generator of the bracket,
\beq
&&(-1)^i X_1\we \dots\we X_{i-1}\we[X_i,V_i]\nn\\
&=&(-1)^i X_1\we \dots\we X_{i-1}\we\Big\{\DE X_i\cdot V_i - \DE V_{i-1} - X_i\we \DE V_i\Big\}\nn\\
&=&(-1)^i \DE X_i\cdot  X_1\we \dots\we X_{i-1}\we X_{i+1}\we \dots \we X_{k+1}\nn\\
&&+(-1)^{i-1} X_1\we \dots\we X_{i-1}\we \DE V_{i-1}\nn\\
&&-(-1)^{i}X_1\we \dots\we X_{i-1}\we X_i\we \DE V_i.
\eeq
By Corollary \ref{al_jU}, 
\beq
(-1)^i \DE X_i\cdot  X_1\we \dots\we X_{i-1}\we X_{i+1}\we \dots \we X_{k+1}
=-\DE X_i\cdot (\al_i\ct U).
\eeq
Thus, we obtain:
\beq
&&\sum_{i=1}^{k}(-1)^i X_1\we \dots\we X_{i-1}\we[X_i,V_i]\nn\\
&=&\sum_{i=1}^{k}-\DE X_i\cdot(\al_i\ct U)\nn\\
&&+\sum_{i=1}^{k}\Big\{(-1)^{i-1} X_1\we \dots\we X_{i-1}\we \DE V_{i-1}
-(-1)^{i}X_1\we \dots\we X_{i-1}\we X_i\we \DE V_i\Big\}\nn\\
&=&\sum_{i=1}^{k}-\DE X_i\cdot(\al_i\ct U)+\DE V_0\nn\\
&&+\sum_{i=1}^{k-1} (-1)^i X_1\we \dots\we X_i\we \DE V_i
-\sum_{i=1}^{k-1} (-1)^i X_1\we \dots\we X_i\we \DE V_i\nn\\
&&-(-1)^k \DE V_k \cdot X_1\we\dots\we X_k\nn\\
&=&\DE U-\sum_{i=1}^{k+1}\DE X_i\cdot(\al_i\ct U),
\eeq
as $V_{k}=X_{k+1}$,\,\, $(-1)^k X_1\we\dots\we X_k=\al_{k+1}\ct U$, and 
\beq
V_0&=&X_1\we\dots\we X_{k+1}\nn\\
&=&U.\nn
\eeq
Therefore, 
\beq
&&\om\ct\sum_{i=1}^{k}\sum_{i<j}(-1)^{i+j}[X_i, X_j]\we 
X_1\we \dots\hat{X_i}\dots\hat{X_j}\dots\we X_{k+1}\nn\\
&=&\om\ct\Big\{\DE U-\sum_{i=1}^{k+1}\DE X_i\cdot\al_i\ct U\Big\}\nn\\
&=&\om\ct \DE U - \sum_{i=1}^{k+1}\DE X_i\cdot \om\ct(\al_i\ct U).
\eeq
Finally, we obtain:
\beq
d\om\ct U &=& -(-1)^{|w|}\DE (\om\ct U)+(-1)^{|w|}\sum_{i=1}^{k+1}\DE X_i\cdot(\al_i\ct(\om\ct U))\nn\\
&&+\om\ct \DE U - \sum_{i=1}^{k+1}\DE X_i\cdot \om\ct(\al_i\ct U)\nn\\
&=&\om\ct \DE U -(-1)^{|w|}\DE (\om\ct U)\nn\\
&&+(-1)^{|w|}\sum_{i=1}^{k+1}\DE X_i\cdot(\al_i\ct(\om\ct U))-
(-1)^{|w|}\sum_{i=1}^{k+1}\DE X_i\cdot(\al_i\ct(\om\ct U))\nn\\
&=&\om\ct \DE U -(-1)^{|w|}\DE (\om\ct U).\nn
\eeq
\\
In the second step we generalize the result to the case where $|w|+1<u$.

%%%%%%%%%%%%%%%%%%%%%%%%%%%%%%%%%%%%%%%%%%%%%%
%
%    GENERAL CASE of
% dw_|U=w_|DU+(-1)^w D(w_|U)
%
%

For any $\theta\in\Ga(\we^{u-1-|\om|}A^*)$, form $\om\we\theta$ is a $u-1$-form, 
therefore we can apply the result proved in the first step:
\beq
d(\om\we\theta)\ct U=(\om\we\theta)\ct \DE U - (-1)^{|\om|+|\theta|}\DE ((\om\we\theta)\ct U).
\label{prop:dw general::main_equ}
\eeq
For the left hand side we have
\beq
d(\om\we\theta)\ct U&=&(d\om\we\theta)\ct U+(-1)^{|\om|}(\om\we d\theta)\ct U\nn\\
&=&\theta\ct(d\om\ct U)+(-1)^{|\om|}d\theta \ct(\om\ct U).\nn
\eeq
The right hand side of (\ref{prop:dw general::main_equ}) can be written as:
\beq
&&(\om\we\theta)\ct \DE U - (-1)^{|\om|+|\theta|}\DE((\om\we\theta)\ct U)\nn\\
&&=\theta\ct(\om\ct \DE U) - (-1)^{|\om|}(-1)^{|\theta|}\DE(\theta\ct(\om\ct U)).\nn
\eeq
Now, $\om\ct U\in\Ga(\we^{|\theta|+1}A)$, so
\beq
(-1)^{|\theta|}\DE (\theta\ct(\om\ct U))=\theta\ct \DE(\om\ct U)-d\theta\ct(\om\ct U).\nn
\eeq
Therefore the relation (\ref{prop:dw general::main_equ}) becomes
\beq
&&\theta\ct(d\om\ct U)+(-1)^{|\om|}d\theta \ct(\om\ct U)\nn\\
&&
=\theta\ct(\om\ct \DE U) + (-1)^{|\om|}d\theta\ct(\om\ct U)-(-1)^{|\om|}\theta\ct \DE(\om\ct U).\nn
\eeq
This gives
\beq
\theta\ct(d\om\ct U)=\theta\ct(\om\ct \DE U)-(-1)^{|\om|}\theta\ct \DE(\om\ct U),\nn
\eeq
which is equivalent to
\beq
\theta\ct\Big\{d\om\ct U-\om\ct\DE U + (-1)^{|\om|}\DE(\om\ct U)\Big\}=0.\label{theta_|all}
\eeq
Since (\ref{theta_|all}) holds for all $\theta\in\Ga(\we^{u-1-|w|}A)$, we must have
\beqs
d\om\ct U-\om\ct\DE U + (-1)^{|\om|} \DE(\om\ct U)=0.\nn
\eeqs
Therefore
\beq
d\om\ct U=\om\ct\DE U - (-1)^{|\om|} \DE(\om\ct U).\label{prop:dw general::main_equ_2}
\eeq
\\
In the third and last step we generalize the result to the case where $U$ may vanish at some points.
For this it is sufficient to show that (\ref{prop:dw general::main_equ_2}) holds 
for $fU, f\in\Cinf(M)$, where $U$ is nowhere vanishing, i.e. that we have
\beq
d\om\ct (fU)=\om\ct\DE (fU) - (-1)^{|\om|} \DE(\om\ct fU).\label{prop:dw general::main_equ_3}
\eeq
We first observe that (\ref{prop:dw general::main_equ_2}) holds for 
$f\in\Ga(\we^0A)=\Cinf(M)$, so
\beq
df\ct U &=& f\ct \DE U - (-1)^{|f|} \DE(f\ct U)\nn\\
&=&f\DE U - \DE(f U),\nn
\eeq
and
\beqs
\om\ct\DE (fU)&=&\om\ct(f\DE U)-\om\ct(df\ct U)\\
&=&f\om\ct\DE U - (df\we\om)\ct U.
\eeqs
It follows from (\ref{prop:dw general::main_equ_2}) that 
\beqs
d(f\om)\ct U &=& (f\om)\ct \DE U - (-1)^{|\om|} \DE((f\om)\ct U).
\eeqs
As $\DE(\om\ct fU)=\DE((f\om)\ct U)$, we can write
\beqs
&&\om\ct \DE(fU) - (-1)^{|\om|} \DE(\om\ct fU) \\
&=&\Big(f\om\ct\DE U - (df\we\om)\ct U\Big)+\Big(d(f\om)\ct U - (f\om)\ct \DE U\Big)\\
&=&d(f\om)\ct U - (df\we\om)\ct U\\
&=&fd\om\ct U\\
&=&d\om\ct(fU).
\eeqs
This finishes the proof.
\end{proof}

%%%%%%%%%%%%%%%%%%%%%%%%%%%%%%%%%%%
%
% Hom/Cohom pairing.
%
%

\subsubsection{Proof of the Theorem \ref{Hom/Cohom pairing}}
\begin{proof}% {\it Proof of the Theorem \ref{Hom/Cohom pairing}}.

We want to demonstrate that for all $\om\in\Ga(\we^{n-k} A^*)$, 
%\beq
%D_0(*\om)=-(-1)^{|\om|}*(d\om).
%\eeq
\beq
\DE_0(\om\ct\La)=-(-1)^{|\om|}d\om\ct\La.\label{thm: hom/cohom pair::expr}
\eeq
By Proposition \ref{prop:dw general} we have 
\beq
d\om\ct\La=\om\ct \DE_0\La-(-1)^{|\om|}\DE_0(\om\ct\La).
\eeq
As $\DE_0\La=0$,
\beq
d\om\ct\La=-(-1)^{|\om|}\DE_0(\om\ct\La),\nn
\eeq
which is equivalent to
\beq
\DE_0*=-(-1)^{n-k}*d.
\eeq
\end{proof}
\begin{acknowledgments}
Authors are very grateful to Ping Xu for introducing us to the subject and for his constant support.
We would also like  to thank Camille Laurent, Aissa Wade and
David Iglesias-Ponte for many useful discussions and suggestions, and
Katherine Hurley for proof-reading this article.
\end{acknowledgments}

\end{document}